\documentclass[11pt]{amsart} 
\usepackage{amscd,amssymb,amsmath} 

\def\T{\mathbb T}
\def\Z{\mathbb Z}

\def\C{\mathbb C}
\def\N{\mathbb N}

\def\H{{\mathcal H}}
\def\K{{\mathcal K}}

\def\U{{\mathcal U}}

\def\id{\operatorname{id}} 
\def\Pol{\operatorname{Pol}} 
\def\Ad{\operatorname{Ad}}

\newtheorem{theo}{Theorem}

\newtheorem{lemm}[theo]{Lemma}

\newtheorem{rema}[theo]{Remark}

\title[Cocycle for quantum $SU(2)$]{A cocycle for the 
comultiplication \\ on the quantum $\mathbf{SU(2)}$ group} 

\author[Szyma\'{n}ski]{Wojciech Szyma\'{n}ski}
\address{Mathematics \\ The University of Newcastle \\
NSW 2308 \\ Australia}
\email{wojciech@frey.newcastle.edu.au} 

\thanks{This work was partially supported by  
Max-Planck-Institute f\"ur Mathematik, Bonn, 
by Universit\"at M\"unster, and by Research 
Grants Committee of the University of Newcastle.}

\date{6 November, 2002} 

\begin{document} 

\begin{abstract} 
For the family $\Delta_q$, $q\in[0,1)$, of the 
$SU_q(2)$-comultiplications on the $C^*$-algebra 
$A\cong C(SU_q(2))$, we show that there exist unitary 
operators $\Omega_q$ such that $\Delta_q(x)=
\Omega_q\Delta_0(x)\Omega_q^*$ for $x\in A$. 
\end{abstract} 

\maketitle 

One common way of producing Hopf $*$-algebras is to twist the 
comultiplication of an already known Hopf $*$-algebra by 
a unitary 2-cocycle or 2-pseudo-cocycle (for example, see \cite{ev}). 
In this note, we show that an analogoues cocycle relation 
exists between the comultiplications on the quantum $SU_q(2)$ 
group corresponding to different values of the parameter 
$q\in[0,1)$. However, our starting point is $q=0$, at which 
$SU_q(2)$ is not a quantum group but only a quantum 
semigroup (with counit). We show that there exist unitary 
operators $\Omega_q$ such that $\Delta_q=
\Ad(\Omega_q)\Delta_0$. Thus, instead of twisting 
existing quantum group structure, $\Omega_q$ `regularizes' 
the comultiplication of a quantum semigroup to a comultiplication
of a quantum group. 

We begin by recalling the definition of the quantum $SU(2)$ group \cite{w1}. 
For $q\in[0,1]$, the $C^*$-algebra $C(SU_q(2))$ of continuous 
functions on the quantum $SU(2)$ group is defined as the 
universal $C^*$-algebra generated by two elements $a$ and $b$, subject 
to the following relations. 
\begin{align*} 
a^*a+b^*b & = I, & ab & =qba, \\  
aa^*+q^2b^*b & = I, & ab^* & =qb^*a, \\   
& & b^*b & =bb^*.   
\end{align*}   
If $q=1$ then this algebra is abelian and isomorphic with 
the $C^*$-algebra of continuous functions on the classical 
group $SU(2)$. If $q\in[0,1)$ then the algebras $C(SU_q(2))$ are 
no longer abelian. However, they are all isomorphic as $C^*$-algebras 
\cite[Theorem A2.2]{w1} to the universal $C^*$-algebra generated by 
two operators $T$ and $S$, subject to the relations 
\begin{align*} 
T^*T & =I, \\ 
S^*S & = SS^*, \\ 
TT^*+S^*S & = I.  
\end{align*} 
We denote this $C^*$-algebra by $A$. By \cite[Theorem A2.2]{w1} and 
\cite[Formulae (A2.2) and (A2.3)]{w1}, there exists an isomorphism 
$\phi_q:C(SU_q(2))\rightarrow A$ given on the generators as follows. 
If $q\in(0,1)$ then 
\begin{align*} 
\phi_q(a) & = \sum_{n=0}^{+\infty}{\Big(}\sqrt{1-q^{2(n+1)}}-
    \sqrt{1-q^{2n}}{\Big)}T^n(T^*)^{n+1}, \\ 
\phi_q(b) & = \sum_{n=0}^{+\infty}q^nT^nS(T^*)^n.  
\end{align*} 
If $q=0$ then  
\begin{align*} 
\phi_0(a) & = T^*,\\ 
\phi_0(b) & = S. 
\end{align*} 
Here and in what follows we agree to write  
$$ T^0=(T^*)^0=I, \; S^0=(S^*)^0=S^*S, \; \text{and} \; S^{-1}=S^*. $$ 

We denote by $J$ the closed 2-sided ideal of $A$ generated by $S$. We have 
\begin{equation}\label{spanofJ} 
J=\overline{\text{span}}\{T^m S^j T^{-n}:m,n\in\N,\,j\in\Z\}. 
\end{equation}  
Thus $J$ is isomorphic to $C(\T)\otimes\K$,   
with $\K$ the compact operators on a separable Hilbert space. The 
quotient is generated by $T+J$ and isomorphic to 
$C(\T)$, the $C^*$-algebra of continuous functions on the circle. 
Thus there is an exact sequence 
$$ 0 \longrightarrow C(\T)\otimes\K\longrightarrow A  
   \longrightarrow C(\T)\longrightarrow 0. $$ 

The ideal $J$ has a natural faithful representation 
on the Hilbert space $\ell^2(\N\times\Z)$, which extends to a 
representation $\rho$ of $A$ constructed by Woronowicz in \cite{w1}. 
Since the ideal $J$ is essential, the representation $\rho$ of $A$ 
is faithful. If $\{\xi(k,m)\;|\;k\in\N,\:m\in\Z\}$ is an orthonormal 
basis of $\ell^2(\N\times\Z)$ then we have 
\begin{align*} 
\rho(T):\xi(k,m) & \mapsto  \xi(k+1,m), \\ 
\rho(S):\xi(k,m) & \mapsto  \delta_{k,0}\xi(k,m+1), \\   
\rho(T^kS(T^*)^k):\xi(k,m) & \mapsto  \xi(k,m+1), \\   
\rho(\phi_q(a)):\xi(k,m) & \mapsto \sqrt{1-q^{2k}}\xi(k-1,m), \\ 
\rho(\phi_q(b)):\xi(k,m) & \mapsto  q^k\xi(k,m+1),  
\end{align*} 
with $\delta$ being the Kronecker symbol, and the convention 
that $\xi(k,m)=0$ whenever $k<0$. 

The irreducible representations of $A$ consist of two families 
$\{\omega_t\;|\;t\in\T\}$ and $\{\rho_t\;|\;t\in\T\}$. Each $\omega_t$ is 
one dimensional and it is determined by $\omega_t(T)=t$, 
$\omega_t(S)=0$. Each $\rho_t$ acts on the infinite dimensional 
Hilbert space $\ell^2(\N)$ with an orthonormal basis $\{\xi(n)\;|\;
n\in\N\}$, and it is determined by $\rho_t(S)\xi(n)=\delta_{n,0}
t\xi(0)$, $\rho_t(T)\xi(n)=\xi(n+1)$. 

We recall the quantum group structure of $SU_q(2)$. The comultiplication 
$$ \Delta:C(SU_q(2))\rightarrow 
   C(SU_q(2))\otimes C(SU_q(2)) $$ 
is a coassociative, unital $C^*$-algebra homomorphism such that 
\begin{align*} 
\Delta(a) & = a\otimes a-qb^*\otimes b, \\ 
\Delta(b) & = b\otimes a+a^*\otimes b.  
\end{align*} 
We transport $\Delta$ to $A$ via the isomorphism 
$\phi_q$, and denote by 
$$ \Delta_q=(\phi_q\otimes\phi_q)\Delta\phi_q^{-1} $$ 
the resulting comultiplication on $A$. If $q=0$ then we have 
\begin{align*} 
\Delta_0(T) & = T\otimes T, \\ 
\Delta_0(S) & = S\otimes T^*+T\otimes S.  
\end{align*} 
For $q>0$ it is a non-trivial task to express $\Delta_q(T)$ 
and $\Delta_q(S)$ in terms of the generators $T,S$. 

\begin{lemm}\label{com-cont} 
For each $x\in A$ the map $q\mapsto\Delta_q(x)$ is norm continuous. 
\end{lemm} 
\begin{proof} 
Since the homomorphism $\Delta_q$ is isometric and the $C^*$-algebra 
$A$ is generated by $S$ and $T$ it sufficies to show that the maps 
$q\mapsto\Delta_q(S)$ and $q\mapsto\Delta_q(T)$ are norm continuous. 

Clearly, the maps $q\mapsto\phi_q(a)$ and $q\mapsto\phi_q(b)$ are 
continuous. It follows that the maps $q\mapsto\Delta_q(\phi_q(a))$ and 
$q\mapsto\Delta_q(\phi_q(b))$ are continuous. Since $1$ is an isolated 
point of the spectrum of $\phi_q(bb^*)$ the spectral projection 
$E_{\phi_q(bb^*)}(\{1\})$ is in $A$. We have $S=
\phi_q(b)E_{\phi_q(bb^*)}(\{1\})$ and hence $\Delta_q(S)=
\Delta_q(\phi_q(b))E_{\Delta_q(\phi_q(bb^*))}(\{1\})$. It follows that 
the map $q\mapsto\Delta_q(S)$ is continuous. For each 
$q\in[0,1)$ the element $|\phi_q(a^*)|$ is invertible and we have 
$T=\phi_q(a^*)|\phi_q(a^*)|^{-1}$. Hence 
$\Delta_q(T)=\Delta_q(\phi_q(a^*))
|\Delta_q(\phi_q(a^*))|^{-1}$ and, consequently, the map 
$q\mapsto\Delta_q(T)$ is continuous. 
\end{proof} 

The counit 
$$ \varepsilon:C(SU_q(2))\rightarrow\C $$ 
is a character such that $\varepsilon(a)=1$ and $\varepsilon(b)=0$. 
Slightly abusing notation we write $\varepsilon=
\varepsilon\circ\phi_q^{-1}$ for the counit on $A$. 
On the generators we have  
$$ \varepsilon(T)=1 \;\;\; \text{and} \;\;\; \varepsilon(S)=0. $$ 
In particular, $\varepsilon$ does not depend on $q$. For all 
$q\in[0,1)$ the following identity, characterizing the counit, holds 
$$ (\varepsilon\otimes\id)\Delta_q=\id=(\id\otimes\varepsilon)\Delta_q. $$ 

We denote by $\Pol_q$ the $*$-subalgebra of $A$ generated by 
$\phi_q(a)$ and $\phi_q(b)$. As shown in \cite{w1,w2}, for each 
$q\in(0,1)$ there exists the coinverse (or antipode)  
$$ \kappa_q:\Pol_q\rightarrow\Pol_q. $$ 
$\Pol_q$ is dense in $A$, while $\Pol_q\cap J$ is dense 
in $J$ (since the ideal $J$ is generated by $\phi_q(b)$). 
We have $\kappa_q(\Pol_q\cap J)=\Pol_q\cap J$. 

Let $\pi:A\rightarrow A/J$ be the canonical surjection and let 
$D=\ker(\pi\otimes\pi)$. It is not difficult to see that $D$ 
coincides with the closure of $A\otimes J+J\otimes A$. 
We define $B$ to be the $C^*$-subalgebra of $A\otimes 
A$ generated by $D$ and $T\otimes T$. 
There exists an exact sequence 
$$ 0 \longrightarrow D \longrightarrow B \longrightarrow C(\T) 
   \longrightarrow 0. $$ 
For all $q\in[0,1)$ we have 
$$ \Delta_q(A)\subseteq B \;\;\;\; \text{ and } 
   \;\;\;\; \Delta_q(J)\subseteq D. $$ 
Indeed, $\pi(\phi_q(a))=\pi(T^*)$ and $\pi(\phi_q(b))=0$  
yield $(\pi\otimes\pi)(\Delta_q(\phi_q(a)))=(\pi\otimes\pi)
(T^*\otimes T^*)$ and $(\pi\otimes\pi)(\Delta_q
(\phi_q(b)))=0$. Consequently, both $\Delta_q(\phi_q(a))$ and 
$\Delta_q(\phi_q(b))$ belong to $B$ and hence $\Delta_q(A)\subseteq B$. 
Since the ideal $J$ of $A$ is generated by $\phi_q(b)$ and 
$\Delta_q(\phi_q(b))\in D$, we have $\Delta_q(J)\subseteq D$. 

Also for all $q\in[0,1)$ we see that 
$$ \Delta_q(x)-\Delta_0(x) \text{ belongs to $D$ for all $x\in A$}. $$ 
Indeed, this is easily checked on the generators $\phi_q(a)$ and 
$\phi_q(b)$ of a dense $*$-subalgebra of $A$. 

\begin{lemm}\label{com-of-J} 
$\overline{D\Delta_q(J)}=D$ for all $q\in[0,1)$.  
\end{lemm} 
\begin{proof} 
At first we observe that $\overline{D\Delta_q(J)}$ 
is a two-sided ideal of $A\otimes A$. We consider separately the 
cases $q\neq0$ and $q=0$. In the former, since $\kappa_q$ 
is an antihomomorphism this follows from the 
fact that $(\kappa_q\otimes\kappa_q)(D\cap(\Pol_q\otimes\Pol_q))=
D\cap(\Pol_q\otimes\Pol_q)$ and $\kappa_q(J\cap\Pol_q)=J\cap\Pol_q$ 
are dense in $D$ and $J$, respectively. In the latter, we notice  
that $\overline{D\Delta_0(J)}$ is the left ideal of $A\otimes A$ 
generated by $\{\Delta_0(ST^{-n}):n=0,1,\ldots\}$. Thus it sufficies 
to check that $\Delta_0(ST^{-n})x\in D\Delta_0(J)$ for $n=0,1,\ldots$ 
and for all $x$ in the set of generators for the algebra $A\otimes A$. 
This can be done directly. For example, 
$$ \Delta_0(ST^{-n})S\otimes T=\delta_{n,0}S^2\otimes I=\delta_{n,0}
   (S\otimes T+S^2T^*\otimes S^*)\Delta_0(S). $$ 
Hence $\overline{D\Delta_q(J)}$ is a two-sided ideal of $A\otimes A$ 
for all $q\in[0,1)$. Now if $\overline{D\Delta_q(J)}\neq D$ then there 
existed an irreducible representation of $D$ vanishing on 
$\overline{D\Delta_q(J)}$. Such a representation of $D$ could be 
extended to an irreducible representation of $A\otimes A$ unitarily equivalent 
to $\omega_t\otimes\rho_z$, $\rho_t\otimes\omega_z$ or $\rho_t\otimes\rho_z$, 
for some $t,z\in\T$. This, however, leads to a contradiction, as one 
can check directly that none of these representations vanishes on 
$D\Delta_q(\phi_q(b))$. 
\end{proof} 

For $q\in[0,1)$ we have a faithful representation $(\rho\otimes
\rho)\Delta_q$ of $A$ on $\ell^2(\N\times\Z)\otimes
\ell^2(\N\times\Z)$. Our first goal is to show that all 
of these representations are unitarily equivalent. To this end 
we need some preparation. 

At first we observe that 
the restriction of $(\rho\otimes\rho)\Delta_q$ to the ideal 
$J$ is still a nondegenerate representation. Indeed, 
otherwise there existed a subspace $\H_0\neq\{0\}$ of 
$\ell^2(\N\times\Z)\otimes\ell^2(\N\times\Z)$ on which 
$(\rho\otimes\rho)\Delta_q(J)$ were zero. If $q\neq0$ then this is 
impossible since one can verify directly that the kernel of 
$(\rho\otimes\rho)\Delta_q(\phi_q(b))$ is trivial. In the case 
$q=0$ this follows from Lemma \ref{com-of-J}. 

Now for $n\in\N$ we define by induction on $k$ numbers 
$\lambda_q(n,k)\in(0,+\infty)$ as 
\begin{align*} 
\lambda_q(n,0) & = 1, \\ 
\lambda_q(n,k+1) & = \lambda_q(n,k)q^{n+2k+1}{\big(}1-q^{2(n+k+1)}
{\big)}^{-1/2}{\big(}1-q^{2(k+1)}{\big)}^{-1/2}. 
\end{align*} 
Since 
$$ {\big(}1-q^{2(n+k+1)}{\big)}^{-1/2}{\big(}1-q^{2(k+1)}{\big)}^{-1/2}
   \leq{\big(}1-q^2{\big)}^{-1} $$ 
it follows that 
\begin{equation}\label{inequality} 
\lambda_q(n,k)\leq{\big(}q^{n+k}\left(1-q^2\right)^{-1}{\big)}^k. 
\end{equation} 
In particular, the series $\displaystyle{\sum_{k=0}^{+\infty}}
\lambda_q(n,k)^2$ converges. We set 
$$ \Lambda_q(n)=\left(\sum_{k=0}^{+\infty}\lambda_q(n,k)^2\right)^{1/2}. $$ 
Now for $n,m\in\N$ such that either $n$ or $m$ (or both) equals $0$ 
and for $i,j\in\Z$ we define a vector $f^q_{n,i,m,j}$ in 
$\ell^2(\N\times\Z)\otimes\ell^2(\N\times\Z)$ as 
$$ f^q_{n,i,m,j}=\sum_{k=0}^{+\infty}\Lambda_q(n+m)^{-1}
   \lambda_q(n+m,k)\xi(n+k,i-k)\otimes\xi(m+k,j+k). $$ 
If $q=0$ then simply $f^0_{n,i,m,j}=\xi(n,i)\otimes\xi(m,j)$. On the 
other hand, if $q\in(0,1)$ then for each vector $\xi(r,x)\otimes\xi(t,y)$, 
$r,t\in\N$ and $x,y\in\Z$, there exists exactly one vector $f^q_{n,i,m,j}$ 
so that $\langle\xi(r,x)\otimes\xi(t,y),\;f^q_{n,i,m,j}\rangle\neq0$. 

\begin{lemm}\label{uq} 
For each $q\in[0,1)$ there exists a unique unitary operator $U_q$ on 
$\ell^2(\N\times\Z)\otimes\ell^2(\N\times\Z)$ such that 
$$ (\rho\otimes\rho)\Delta_q(x)=U_q((\rho\otimes\rho)\Delta_0(x))U_q^* $$ 
for all $x\in A$, and 
$$ U_q:f^0_{n,i,m,j}\mapsto f^q_{n,i,m,j} $$ 
for all $i,j\in\Z$ and $n,m\in\N$ with $nm=0$. 
\end{lemm} 
\begin{proof} 
For each $q\in[0,1)$, $0$ is an isolated point of the spectrum 
of $\phi_q(a^*a)$ and the corresponding spectral projection is 
$S^*S$. Thus, $(\rho\otimes\rho)\Delta_q(S^*S)$ is the 
projection onto the kernel of $(\rho\otimes\rho)\Delta_q(\phi_q(a))$. 
It is not difficult to verify that this kernel consists of all 
vectors $\sum c(n,i,m,j)\xi(n,i)\otimes\xi(m,j)$ in $\ell^2(\N\times\Z)
\otimes\ell^2(\N\times\Z)$ such that 
$$ c(n+1,i-1,m+1,j+1)=q^{n+m+1}{\big(}1-q^{2(n+1)}{\big)}^{-1/2}
   {\big(}1-q^{2(m+1)}{\big)}^{-1/2}c(n,i,m,j). $$ 
This identity implies that 
vectors $\{f^q_{n,i,m,j}\;|\;n,m\in\N,\:nm=0,\:i,j\in\Z\}$ form an 
orthonormal basis of $(\rho\otimes\rho)\Delta_q(S^*S)$. 

Since $\Delta_q(S)=\Delta_q(\phi_q(b))\Delta_q(S^*S)$ 
we find that 
$$ (\rho\otimes\rho)\Delta_q(S)f^q_{n,i,0,j}=\Lambda_q(n)^{-1}
   \sum_{k=0}^{+\infty}\xi((n+1)+k,i-k)\otimes\xi(k,(j+1)+k) $$  
$$ {\Big(}\lambda_q(n,k+1)q^{n+k+1}{\big(}1-q^{2(k+1)}
   {\big)}^{1/2}+\lambda_q(n,k)q^k{\big(}1-q^{2(n+k+1)}{\big)}^{1/2}{\Big)}. $$ 
This implies that $\langle(\rho\otimes\rho)\Delta_q(S)
f^q_{n,i,0,j},\;f^q_{r,x,t,y}\rangle=0$ if $(r,x,t,y)
\neq(n+1,i,0,j+1)$. Thus, since $(\rho\otimes\rho)\Delta_q
(S)$ is a partial unitary with domain and range $(\rho\otimes\rho)
\Delta_q(S^*S)$, we see that $(\rho\otimes\rho)\Delta_q(S)
f^q_{n,i,0,j}$ is a scalar multiple of $f^q_{n+1,i,0,j+1}$. 
However, it is clear from the series expansions of 
$f^q_{n+1,i,0,j+1}$ and $(\rho\otimes\rho)\Delta_q(S)
f^q_{n,i,0,j}$ that this scalar is a positive real number and, 
hence, equal to $1$. A similar reasoning yields the value of 
$(\rho\otimes\rho)\Delta_q(S)f^q_{0,i,m,j}$ 
for $m\neq0$. We have 
\begin{align} 
(\rho\otimes\rho)\Delta_q(S)f^q_{0,i,m,j} & = f^q_{0,i+1,m-1,j}, \nonumber \\ 
(\rho\otimes\rho)\Delta_q(S)f^q_{n,i,0,j} & = f^q_{n+1,i,0,j+1}, \nonumber 
\end{align} 
for $i,j\in\Z$ and $n,m\in\N$, $m\neq0$. 
We now define an isometry $\widetilde{U}_q$ with domain $(\rho\otimes\rho)
\Delta_0(S^*S)$ and range $(\rho\otimes\rho)\Delta_q(S^*S)$ by 
$$ \widetilde{U}_q:f^0_{n,i,m,j}\mapsto f^q_{n,i,m,j}. $$ 
By construction, we have 
$$ (\rho\otimes\rho)\Delta_q(S)=\widetilde{U}_q
   ((\rho\otimes\rho)\Delta_0(S))\widetilde{U}_q^*. $$ 
Since the restriction of $(\rho\otimes\rho)\Delta_q$ to  
$J$ is still a nondegenerate representation we can  
extend $\widetilde{U}_q$ to a unitary operator $U_q$ on 
$\ell^2(\N\times\Z)\otimes\ell^2(\N\times\Z)$ by setting 
\begin{equation}\label{uqop} 
U_q=\sum_{k=0}^{+\infty}((\rho\otimes\rho)\Delta_q
(T)^k)\widetilde{U}_q(\rho\otimes\rho)
\Delta_0(T^*)^k.  
\end{equation} 
It follows immediatelly from the construction and the property 
of $\widetilde{U}_q$ that 
\begin{eqnarray*} 
(\rho\otimes\rho)\Delta_q(T) & = & 
U_q((\rho\otimes\rho)\Delta_0(T))U_q^*, \\ 
(\rho\otimes\rho)\Delta_q(S) & = & 
U_q((\rho\otimes\rho)\Delta_0(S))U_q^*.   
\end{eqnarray*} 
Furthermore, (\ref{uqop}) is the unique extension of $\widetilde{U}_q$ 
which satisfies the above two identities. 
Since the $C^*$-algebra $A$ is generated by 
$T$ and $S$ it follows that 
$$ (\rho\otimes\rho)\Delta_q(x)=U_q((\rho\otimes\rho)\Delta_0(x))U_q^* $$ 
for all $x\in A$, as desired. 
\end{proof} 

\begin{lemm}\label{uqmult}  
For each $q\in[0,1)$, the unitary $U_q$ of (\ref{uqop}) satisfies 
\begin{description} 
\item[(i)] $U_q(\rho\otimes\rho)(D)\subseteq(\rho\otimes\rho)
(D)$ and $(\rho\otimes\rho)(D)U_q\subseteq(\rho\otimes\rho)(D)$. 
\item[(ii)] $U_q(\rho\otimes\rho)(J\otimes J)\subseteq(\rho\otimes\rho)
(J\otimes J)$ and $(\rho\otimes\rho)(J\otimes J)U_q\subseteq
(\rho\otimes\rho)(J\otimes J)$. 
\item[(iii)] For each $w\in D$ the maps $q\mapsto 
U_q(\rho\otimes\rho)(w)$ and $q\mapsto(\rho\otimes\rho)(w)U_q$ 
are norm continuous. 
\end{description}  
Furthermore, the partial isometry $\widetilde{U}_q:f^0_{n,i,m,j}
\mapsto f^q_{n,i,m,j}$ belongs to $(\rho\otimes\rho)(D)$. 
\end{lemm} 
\begin{proof} 
An elementary calculation on the basis shows that $\widetilde{U}_q$ is 
equal to the sum of the series  
\begin{equation}\label{uqtilde} 
(\rho\otimes\rho){\Big(}I\otimes I-TT^*\otimes TT^*    
   +\sum_{\stackrel{n,m\in\N}{nm=0}}\left(\Lambda_q(n+m)^{-1}-1\right)
   T^n S^*ST^{-n}\otimes T^m S^*ST^{-m}   
\end{equation} 
$$ +\sum_{\stackrel{n,m\in\N}{nm=0}}\Lambda_q(n+m)^{-1}
   \sum_{k=1}^{+\infty}\lambda_q(n+m,k)T^{n+k}S^{-k}T^{-n}\otimes T^{m+k}
   S^{k}T^{-m}{\Big)}, $$ 
and hence belongs to $(\rho\otimes\rho)(D)$. Inequality (\ref{inequality}) 
implies that this series converges in norm and, furthermore, that the map 
$q\mapsto\widetilde{U}_q$ is norm continuous.  

By virtue of identity (\ref{spanofJ}) and Lemma \ref{com-of-J}, to prove 
$U_q(\rho\otimes\rho)(D)\subseteq(\rho\otimes\rho)(D)$ 
it sufficies to show that $U_q(\rho\otimes\rho)(\Delta_0(T^mS^iT^{-n})
w)$ belongs to $(\rho\otimes\rho)(D)$ for all 
$w\in D$, $m,n\in\N$, $i\in\Z$. Indeed, formula (\ref{uqop}) implies that 
$$ U_q(\rho\otimes\rho)(\Delta_0(T^mS^iT^{-n})w)=(\rho\otimes\rho)
   (\Delta_q(T^m))\widetilde{U}_q(\rho\otimes\rho)(\Delta_0(T^{-m}))
   (\rho\otimes\rho)(\Delta_0(T^mS^iT^{-n})w). $$ 
This element belongs to $(\rho\otimes\rho)(D)$, since $\widetilde{U}_q$ 
is in $(\rho\otimes\rho)(D)$. Furthermore, the map $q\mapsto U_q(\rho\otimes
\rho)(v)$ is norm continuous for $v=\Delta_0(T^mS^iT^{-n})w$, 
and consequently it is norm continuous for all $w\in D$. 
A similar argument yields the inclusion $(\rho\otimes\rho)(D)
U_q\subseteq(\rho\otimes\rho)(D)$ and norm continuity of the map 
$q\mapsto(\rho\otimes\rho)(w)U_q$, $w\in D$. This prooves (i) 
and (iii). Part (ii) follows from (i), since $J\otimes J$ is an 
ideal od $D$. 
\end{proof} 

We denote by $M(J\otimes J)$ the multiplier $C^*$-algebra of $J\otimes J$ 
and by $\U(M(J\otimes J))$ the unitary group of $M(J\otimes J)$. Since 
$J\otimes J$ is an essential ideal of $A\otimes A$ there exists a canonical 
imbedding of $A\otimes A$ into $M(J\otimes J)$ and, hence, we may identify 
$A\otimes A$ with a subalgebra of $M(J\otimes J)$. 

Since $\rho\otimes\rho$ is a faithful, non-degenerate representation of 
$J\otimes J$, it uniquely extends to a faithful representation of 
$M(J\otimes J)$, which we still denote $\rho\otimes\rho$. If $U_q$ is 
the unitary operator defined in (\ref{uqop}) then $U_q$ belongs to 
$(\rho\otimes\rho)(M(J\otimes J))$ by Lemma \ref{uqmult}. 
Thus we can define   
\begin{align} 
\Omega_q &=(\rho\otimes\rho)^{-1}(U_q), \label{omegaq} \\ 
\widetilde{\Omega}_q &=(\rho\otimes\rho)^{-1}(\widetilde{U}_q). 
\label{omegatildeq} 
\end{align}  
It follows that $\Omega_q$ is a unitary in $M(J\otimes J)$ and 
$\widetilde{\Omega}_q$ is a partial isometry in $D$. 

Since $(\varepsilon\otimes\id)(D)=J=(\id\otimes\varepsilon)(D)$, 
there exist unique homomorphism extensions (denoted by the same symbols 
as the original maps): 
\begin{align*} 
\varepsilon\otimes\id &: D\rightarrow J & \text{ to } \;\;\;\; 
   \varepsilon\otimes\id &: M(D)\rightarrow M(J), \\ 
\id\otimes\varepsilon &: D\rightarrow J & \text{ to } \;\;\;\; 
   \id\otimes\varepsilon &: M(D)\rightarrow M(J). 
\end{align*} 
By Lemma \ref{com-of-J}, $(\Delta_q\otimes\id)(J\otimes J)\cdot D\otimes J$ 
is dense in $D\otimes J$, and $(\id\otimes\Delta_q)(J\otimes J)\cdot J\otimes D$ 
is dense in $J\otimes D$. Thus there exist unique homomorphism extensions: 
\begin{align*}  
\Delta_q\otimes\id &: J\otimes J\rightarrow D\otimes J & \text{ to } \;\;\;\;
   \Delta_q\otimes\id &: M(J\otimes J)\rightarrow M(D\otimes J), \\ 
\id\otimes\Delta_q &: J\otimes J\rightarrow J\otimes D & \text{ to } \;\;\;\;
   \id\otimes\Delta_q &: M(J\otimes J)\rightarrow M(J\otimes D). 
\end{align*} 

\begin{theo}\label{cocycle} 
For $q\in[0,1)$ let $\Omega_q\in\U(M(J\otimes J))$ be the unitary 
defined in (\ref{uqop}). Then the following conditions hold true. 
\begin{description} 
\item[(i)] $\Ad\Omega_q$ is an automorphism of $B$. 
\item[(ii)] $\Delta_q(x)=\Omega_q\Delta_0(x)\Omega_q^*$ for all $x\in A$. 
\item[(iii)] The map $[0,1)\ni q\mapsto\Omega_q\in\U(M(J\otimes J))$ is 
continuous in the strict topology.  
\item[(iv)] $(\varepsilon\otimes\id)(\Omega_q)=I=
(\id\otimes\varepsilon)(\Omega_q)$. 
\end{description} 
\end{theo} 
\begin{proof} 
The unitary $\Omega_q$ satisfies (ii) and $\Omega_q D\Omega_q^*=D$ by 
by Lemma \ref{uq}. Since $\Omega_q(T\otimes T)\Omega_q^*=\Delta_q(T)$ 
equals $T\otimes T$ plus an element of $D$, it follows that 
$\Ad\Omega_q(B)=B$. That is, condition (i) is satisfied. 
Condition (iii) follows from Lemma \ref{uqmult}. 

To prove $(\varepsilon\otimes\id)(\Omega_q)=I$ it sufficies to show that 
this equality holds when multiplied on the right by an arbitrary 
element $x=(\varepsilon\otimes\id)(\Delta_0(x))$ of $J$. By 
(\ref{spanofJ}), it is enough to take $x=T^mS^iT^{-n}$. We have 
$$ (\varepsilon\otimes\id)(\Omega_q)(\varepsilon\otimes\id)\left(
   \Delta_0(T^mS^iT^{-n})\right)=(\varepsilon\otimes\id)\left
   (\Omega_q\Delta_0(T^mS^iT^{-n})\right) $$ 
$$ =(\varepsilon\otimes\id)\left(\Delta_q(T^m)\widetilde{\Omega}_q
   \Delta_0(T^{-m})\Delta_0(T^mS^iT^{-n})\right). $$ 
Using (\ref{uqtilde}) and (\ref{omegatildeq}) one can  
show that this element equals $T^mS^iT^{-n}$. Thus $(\varepsilon\otimes\id)
(\Omega_q)=I$. Similarly, one shows that $I=(\id\otimes\varepsilon)(\Omega_q)$, 
and (iv) is proved. 
\end{proof} 

\begin{rema}
\rm Since $\Delta_q$ is coassociative it easily follows from 
condition (ii) of Theorem \ref{cocycle} that the unitary element 
\begin{equation}
(\id\otimes\Delta_0)(\Omega_q^*)(I\otimes\Omega_q^*)
(\Omega_q\otimes I)(\Delta_0\otimes\id)(\Omega_q) 
\end{equation} 
of $M(J\otimes J\otimes J)$ belongs to the relative commutant 
of $(\Delta_0\otimes\id)\Delta_0(A)$. Thus, following 
\cite{ev}, we may call $\Omega_q$ a 2-pseudo-cocycle for 
$\Delta_0$. We still do not know if $\Omega_q$ is a 
2-cocycle for $\Delta_0$; that is, if it satisfies 
\begin{equation}
(\Omega_q\otimes I)(\Delta_0\otimes\id)(\Omega_q)=
(I\otimes\Omega_q)(\id\otimes\Delta_0)(\Omega_q). 
\end{equation}
\end{rema}

\end{document}